\documentclass[12pt]{elsarticle}

\usepackage{amssymb}
 \usepackage{amsthm}

\usepackage[utf8x]{inputenc}
\usepackage[english]{babel}
\usepackage{amsmath}
\usepackage{colonequals}
\usepackage[table,usenames,dvipsnames]{xcolor}
\usepackage{hyperref}
\usepackage{subcaption}
\usepackage{tikz}

\newcommand{\R}{\mathbb{R}}
\newcommand{\N}{\mathbb{N}}

\newcommand{\op}[1]{\operatorname{#1}}
\providecommand{\argmin}[1]{\underset{#1}{\op{arg\,min\,}}}
\newcommand{\half}{{\textstyle \frac{1}{2}}}

\global\long\def\cA{A}
\global\long\def\cP{\mathcal{P}}

\global\long\def\cS{\mathcal{S}}

\global\long\def\cP{\mathcal{P}}
\global\long\def\cO{\mathcal{O}}
\global\long\def\cJ{\mathcal{J}}

\newcommand{\ECH}{\cJ_{\text{CH}}}

\global\long\def\cE{\mathcal{E}}
\global\long\def\cM{\mathcal{M}}
\global\long\def\cP{\mathcal{P}}

\newcommand{\q}{\textbf{q}}            
\newcommand{\vel}{\textbf{p}}          
\newcommand{\cD}{\mathcal D}           
\newcommand{\TD}{\mathrm{T}\cD}        
\newcommand{\EXP}{\mathbb{E}}          
\newcommand{\bgamma}{\mathbf r}        
\newcommand{\spn}{\operatorname{span}}  

\newcommand{\ce}{\mathbf e}
\newcommand{\Te}{\theta}
\newcommand{\fric}{\gamma}
\newcommand{\qA}{{\mathbf Q}}
\newcommand{\qG}{{\mathbf G}}
\newcommand{\qB}{{\mathbf B}}
\newcommand{\sigmaLangevin}{\widehat{\sigma}_{\rm Langevin}}

\usepackage{todonotes}
\usepackage[normalem]{ulem}

\newtheorem{theorem}{Theorem}

\graphicspath{{gfx/}}

\begin{document}

\begin{frontmatter}

\title{Free energy computation of particles with membrane-mediated interactions via Langevin dynamics}

\author{Tobias Kies\fnref{labelKies}}\ead{tobias.kies@fu-berlin.de}
\fntext[labelKies]{IAV GmbH, Nordhoffstr. 5, 38559 Gifhorn, Germany}
\author{Carsten Gr\"{a}ser\fnref{labelRest}}\ead{graeser@mi.fu-berlin.de}
\fntext[labelRest]{Freie Universität Berlin, Institute of Mathematics, Arnimallee 6, 14195 Berlin, Germany}
\author{Luigi Delle Site \fnref{labelRest}}\ead{luigi.dellesite@fu-berlin.de}
\author{Ralf Kornhuber \fnref{labelRest}} \ead{ralf.kornhuber@fu-berlin.de}

\begin{abstract}
We apply well-established concepts of Langevin sampling to derive a new class of algorithms 
for the efficient computation of free energy differences of fluctuating particles 
embedded in a 'fast' membrane, i.e., a membrane  that instantaneously adapts to varying particle positions.
A geometric  potential accounting for membrane-mediated particle interaction is derived 
in the framework of variational hybrid models for particles in membranes.
Recent explicit representations of the gradient of the geometric interaction potential allows to apply  
well-known gradient based Markov Chain Monte-Carlo (MCDC) methods such as Langevin-based sampling.
\end{abstract}




\end{frontmatter}



\section{Introduction} \label{sec:INTRO}

%
The interplay of proteins and curvature of lipid bilayers is well-known to regulate cell
morphology and a variety of cellular functions, such as trafficking or signal detection~\cite{Lip91, McMahonGallupNature2005, SimVot15}.
Microscopic causes, such as hydrophobic mismatch of proteins and amphiphilic lipids, 
may have macroscopic effects, such as budding or fission. 
For example, the membrane remodeling during clathrin-mediated endocytosis 
involves concerted actions of highly specialized membrane proteins 
that can both sense and create membrane curvature, cf., e.g., \cite{haucke2018membrane} 
and the literature cited therein.

%
%
A well-established approach to the modeling of particles in lipid membranes is based on \emph{coarse-grained molecular dynamics}. 
In order to overcome well-known limitations of classical molecular dynamics with respect  to length and time scales, the membrane constituents,
i.e., lipids and proteins, are represented by short chains of beads, cf., e.g., \cite{brannigan2008model,laradji2011coarse,SaundersVoth2014,WD10}.
On the macroscopic side of the model hierarchy, there are  \emph{pure continuum models}
based on the fundamental Canham--Helfrich model of lipid membranes~\cite{Can70,Hel73}
and the representation of proteins by areal concentrations. 
The mutual coupling of particles and membrane is described by concentration-dependent 
mechanical properties of the membrane, as, e.g., bending rigidities or spontaneous curvature,
and line energies of the concentrations associated with phase boundaries~\cite{Lip92,JulLip96}.

%
\emph{Hybrid models} are intended to bridge the gap between coarse-grained molecular dynamics based approaches
that  still have certain limitations in terms of the accessible time and
length scales, and pure continuum models that are unable to incorporate the effect of small particle counts. 
These models are hybrid in the sense that the continuous Canham--Helfrich model of lipid membranes
is coupled to  a finite number of discrete particles represented as finite size or point-like objects.
There is a rich literature on coupling conditions typically prescribing
contour and slope of the membrane either at particle boundaries~\cite{GouBruPin93,WeiKozHel98,SchweitzerShemeshKozlov2017}
or in single points~\cite{KimNeuOst98,DomFou99,DomFou02,MarMis02,BarFou03,WeiDes13}.
In a recently developed variational approach to hybrid models~\cite{ElliottEtAl2016},
see also~\cite{elliott2019second,elliott2019small},
such coupling conditions take the role of constraints in energy minimization.
While hybrid models are often formulated in the zero-temperature limit, 
effects of \emph{thermal fluctuations}, mostly of the membrane, are about to attract more and more
attention~\cite{Netz97,gov2004membrane,naji2009hybrid,duncan2015multiscale,sigurdsson2016hydrodynamic}.

%
In this paper, we aim at \emph{macroscopic properties of fluctuating particles in  membranes}.
To this end, we apply well-established concepts of Langevin sampling of free energy differences cf., e.g., \cite{bussi2007accurate,latorre,lelievre2012langevin,leimkuhler2013robust},
to particles with membrane-mediated interactions.
Assuming that the membrane undergoes only small deformations (Monge gauge) 
and is 'fast' in the sense that its shape instantaneously adapts to varying particle positions,
we suggest a geometric  potential  accounting for the membrane-mediated interaction of  particles.
This geometric interaction potential  can be numerically evaluated by constrained minimization of the associated membrane energy
or, equivalently, by approximate solution of a  fourth-order partial differential equation,  the corresponding Euler--Lagrange equation.
The computation of free energy differences is then performed by 
the classical energy perturbation method due to Zwanzig~\cite{zwanzig1954high}, \cite[Section~1.3]{lelievre}
leading to corresponding  sampling problems for integrals over high-dimensional phase space.
It is well-known, cf., e.g.,  Durmus et al.~\cite{durmus},
that \emph{gradient based Markov Chain Monte-Carlo (MCDC) methods} such as Langevin-based samplers~\cite{bussi2007accurate,latorre,lelievre2012langevin,leimkuhler2013robust}, 
typically have better convergence properties than gradient free  Metropolis-type algorithms.
Exploiting recent results on derivatives of the geometric interaction potential~\cite{GraeserKies2017a}, 
we are finally able to introduce  novel Langevin-based sampling methods 
for the computation of free energy differences of fluctuating particles in 'fast' membranes.

The paper is organized as follows.
First, we give a brief overview on variational hybrid methods for particles in membranes
and formally define  the geometric interaction potential for a general class of coupling conditions.
Then we consider finite-size particles imposing a certain contour and slope on the membrane  at their boundary
as a concrete and practically relevant example,
derive the corresponding geometric interaction potential, and present a computationally feasible representation of its  derivative.
On this background, we formulate an overdamped Langevin equation for an interaction potential accounting 
for both geometric and direct particle--particle interaction
and  briefly recall the classical energy perturbation method for the computation of free energy differences.
Finally, we formulate  our Langevin based sampling method and discuss its main theoretical properties.

\section{Variational hybrid modeling of particles in membranes} \label{sec:VMPM}

\subsection{Geometric interaction potential} \label{subsec:GEOPOT}
%
%
We consider the  fundamental Canham--Helfrich bending energy~\cite{Can70,Hel73}
\[
\ECH(\cS) =   \int_{\cS} \half \kappa(H - c_0)^2 + \kappa_G K \; d \cS,
\]
where the surface $\cS \subset \R^3$ with surface element $d\cS$  is representing the membrane, 
$H$ and $K$ stand for mean and Gaussian curvature, 
$\kappa>0$  and $\kappa_G>0$ are the corresponding bending rigidities,
and $c_0$ is a preferred spontaneous mean curvature.
In the special case $\kappa=1$, $c_0=0$, and $\kappa_G = 0$, we obtain the Willmore energy
which plays an important role in minimal surface theory~\cite{willmore1965note}.

Excluding topological changes, we from now on  ignore Gaussian curvature in light of the Gauss--Bonnet theorem.
In addition, we assume  that $\cS = \{ (x,u(x))\;|\; x \in \Omega\}$ can be parametrized as a graph 
over an open, bounded, non-empty reference domain $\Omega \subset \R^2$ with sufficiently smooth boundary $\partial \Omega$
and that  $\cS$ is almost flat, i.e. $|\nabla u| \ll 1$.
Then it is justified to approximate $\cJ_{CH}(\cS)$ by  the linearized  Canham--Helfrich energy (Monge gauge)
\begin{equation} \label{eq:CH}
\cJ_{\Omega}(u) = \half \int_\Omega   \kappa (\Delta u)^2 +   \sigma |\nabla u|^2\; dx ,
\end{equation}
where the bending energy is additionally supplemented by a surface energy $\int_{\Omega}\sigma |\nabla u|^2\; dx$
that is associated with membrane tension $\sigma\geq 0$.

 We also consider a finite number $N \in \N$ of particles  $B_{\q_i}$  each of 
 which is represented by a position vector $\q_i \in \R^q$, $i=1,\dots,N$.
 We assume that  $\q =(\q_i)_{i=1}^N$ 
 is contained in a given configuration space $\cD\subset \R^{q \times N}$ of feasible particle positions.
In particular, $\cD$ is chosen such that  for all $\q \in \cD$  the membrane equipped with particles with positions $\q$  
can be parametrized over the domain $\Omega_{\q}\subset \Omega$. We assume that $\Omega_{\q}$ is sufficiently regular 
so that the energy $\cJ_{\Omega_\q}(u) $ is well-defined for deformations $u$ taken from the Sobolev space $H^2(\Omega_{\q})$ 
consisting of functions in $L^2(\Omega_{\q})$ that have weak derivatives of second order in $L^2(\Omega_{\q})$.
With $\frac{\partial}{\partial n}v$ denoting the weak derivative of $v \in H^2(\Omega_{\q})$ 
in the direction of the outward normal on $\partial \Omega$,
we assume that $\cD$ is chosen such that   $\cJ_{\Omega_\q}$ is coercive on
 $H_{\q} = \{ v \in H^2(\Omega_{\q})\;|\;v = \frac{\partial}{\partial n} v = 0 \text{ on }\partial \Omega\}$ for $\q \in \cD$.
Notice that such kind of subspaces of $H^2(\Omega_{\q})$ are associated with the classical concept of weak solutions
of partial differential equations as, e.g., 
arising as Euler--Lagrange equations of minimization problems with energies of the form $\cJ_{\Omega_\q}$.

Let the parameters $\bgamma_i\in \R^s$ describe additional degrees of freedom of the particle $B_{\q_i}$
that are automatically determined in course of minimizing the energy of the  membrane for fixed $\q \in \cD$.
As possible examples, one might think of the height or tilt of particles which typically vary freely with membrane deformation.
The coupling of particles with the membrane is then performed by the abstract  conditions
\begin{equation} \label{eq:INTERCOND}
  \exists \bgamma=(\bgamma_i)_{i=1}^N \in \R^{s \times N}: \qquad g_i(u,\bgamma_i; \q_i)=0, \quad i=1,\dots,N.
\end{equation}
Here, $g_i(\cdot; \q_i)$ stands for suitable affine linear continuous mappings from $H_{\q}\times \R^{s}$ to suitable discrete or function spaces.
For a variety of possible selections of the mappings $g_i$, we refer to  \cite{ElliottEtAl2016} 
and also to the example  \eqref{eq:COUPLING} below.
For given particle positions $\q \in \cD$, the space of feasible membrane deformations is therefore given by
\begin{align}\label{eq:constrained_subspace}
  W_{\q} =  \{ v \in H_{\q} \;|\; \exists \bgamma \in \R^{s \times N} g_i(v, \bgamma_i; \q_i)=0 \;\forall  i = 1,\dots, N\} .
\end{align}
Assuming that there is at least one $u \in H_{\q}$
satisfying the constraints~\eqref{eq:INTERCOND},
we find that $W_{\q}$ is a non-empty, affine linear subspace
such that the existence of a unique
\[
u_{\q} = \argmin{ v \in W_{\q}}   \cJ_{\Omega_{\q}}(v)
\]
follows from the Lax--Milgram lemma.
On this background, the  \emph{geometric interaction potential of particles in membranes} is well-defined according to
\[
\cM(\q) = \cJ_{\Omega_{\q}}(u_{\q}) = \underset{v \in  W_{\q} }{\min}\cJ_{\Omega_{\q}}(v),\qquad \q \in \cD.
\]
Though all considerations and algorithms derived below directly extend to this general setting, 
we from now on concentrate on the special case of finite-size particles with curve constraints 
to fix the ideas.
%
%
\subsection{Finite-size particles with curve constraints} \label{subsec:curveconstraints}
Transmembrane proteins are interacting with the membrane curvature by the shape of the hydrophobic belt~\cite{McMahonGallupNature2005}.
Other particles, like FCHo proteins are acting as active or passive scaffolds~\cite{McMahonGallupNature2005} 
or might be partially wrapped due to adhesion energy~\cite{BahramiEtAl2014,KoltoverEtAl1999}.
All these phenomena can be captured by describing the proteins as finite-size, rigid particles, 
imposing a specific contour and slope of the membrane at their boundary.
We will use this example in order to illustrate  the abstract coupling conditions \eqref{eq:INTERCOND}.

Consider $N$ reference particles expressed by  $N$ non-empty, bounded, open  sets $B_i \subset \R^2$ with $0 \in B_i$ 
and sufficiently smooth boundaries $\Gamma_i = \partial B_i$, $i= 1,\dots,N$. 
Moving these reference particles $B_i$ around, we obtain the particles  $B_{\q_i}$,
\[
B_{\q_i} \colonequals \Phi_{\q_i} B_i ,\quad \Phi_{\q_i}(y) = X_i + 
\left( \begin{matrix}
\cos(\alpha_i) & -\sin(\alpha_i)\\
\sin(\alpha_i) & \cos(\alpha_i)
\end{matrix} \right)y, \quad i=1,\dots, N,
\]
with position vectors $\q_i=(X_i,\alpha_i) \in \R^2\times \R$ 
composed of lateral translations $X_i\in \R^2$ and rotation angles $\alpha_i$.
The boundaries of the moved particles $B_{\q_i}$ are denoted by
$\Gamma_{\q_i} = \partial B_{\q_i} = \Phi_{\q_i} \Gamma_i$.
The translation and rotation of the particle is
illustrated in the left picture of Figure~\ref{fig:particle}.

We identify $\R^2\times \R$ with $\R^3$ and introduce the joint position vector $\q = (\q_i)_{i=1}^N\in \R^{3N}$ and
the subset 
\[\Omega_{\q}= \Omega \setminus \bigcup_{i=1}^N \overline{B}_{\q_i} \subset \Omega
\] 
occupied by the membrane. 
We assume that the configuration space  
\[
  \cD=\{ \q \in \R^{3 N}\;|\; \overline{B}_{\q_i}\subset \Omega\quad \forall i, \;
\overline{B}_{\q_i}\cap \overline{B}_{\q_j} = \emptyset \quad  \forall i\neq j\}
\]
of feasible particle positions, 
that are contained in $\Omega$ and do  neither touch each other nor the boundary of $\Omega$, is non-empty. 

We assume that the contour and the slope of the membrane 
along the reference boundaries $\Gamma_i$ is prescribed by given functions $h_i$ and $s_i$, respectively.
Then the contour and slope along the moved particle
boundaries $\Gamma_{\q_i}$ are given by $h_i(\Phi_{\q_i}^{-1}(\cdot))$
and $s_i(\Phi_{\q_i}^{-1}(\cdot))$, respectively.
In addition, particles can freely move up and down and tilt with the ambient membrane. 
This suggests coupling conditions of the form  \eqref{eq:INTERCOND} with $\bgamma_i=(Z_i, \beta_i)$
for $i=1,\dots,N$,
where $Z_i \in \R$ is the height and $\beta_i=(\beta_{i,1}, \beta_{i,2}) \in \R^2$
are linearized tilt angles
around the axes $x_1$, $x_2$ of the particle $B_{\q_i}$.
Combination with prescribed contour and slope leads to  the choice
\begin{equation} \label{eq:COUPLING}
 g_i(u, \bgamma_i; \q_i )(x) = \left( \begin{matrix}
   u(x) - \left(h_i(\Phi_{\q_i}^{-1} (x))+ Z_i + \beta_i \cdot (x-X_i) \right)\\
 {\textstyle \frac{\partial}{\partial n}} u(x) - \left(s_i(\Phi_{\q_i}^{-1} (x)) + \beta_i \cdot n(x) \right) \end{matrix}  \right)
\quad x \in \Gamma_{\q_i}                      
\end{equation}
where
 $x \cdot y$ denotes the Euclidean inner product in $\R^2$
and $n(x)$ stands for the outward normal on $\Omega_{\q}$ in $x \in \Gamma_{\q_i}$.
Note that we made use of the function
\begin{align}\label{eq:height_tilt_function}
  x \mapsto Z_i + \beta_i \cdot (x-X_i)
\end{align} 
and its normal derivative $x \mapsto \beta_i \cdot n(x)$
representing the vertical translation by $Z_i$ and
tilt by $\beta_i$ of the
particle located at $X_i$.
The conditions on contour $h_i$ and slope $s_i$ at a reference boundary $\Gamma_i$ are illustrated
in the right picture of Figure~\ref{fig:particle} in the case of a transmembrane protein.
In practical applications, the functions  $h_i$ and $s_i$ 
can be derived from local molecular dynamics simulations~\cite{NetzZendeHroudLoche20}.

\setlength{\unitlength}{1cm} 
\begin{figure}[ht]
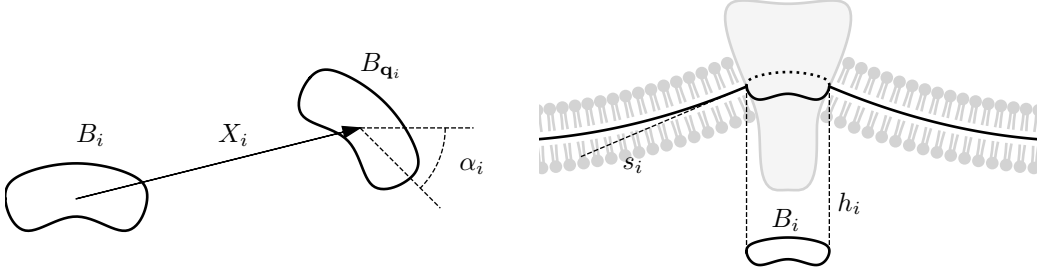

  \centering
  \input{particle_top_view.pgf}%
  \hfill
  \input{particle_side_view.pgf}%

  \caption{%
    Left: Top view of reference particle $B_i$ and moved particle $B_{\q_i}$ with position $\q_i = (X_i,\alpha_i)$.
    Right: Side view of reference particle $B_i$ with contour $h_i$ 
    and slope $s_i$ prescribed at $\Gamma_i$.}.
  \label{fig:particle}
\end{figure}

We now consider the minimization of $\cJ_{\Omega_\q}$ on the
affine subspace $W_{\q} \subset H_{\q}$ defined by~\eqref{eq:constrained_subspace}
for $g_i(\cdot;\cdot)$ given by~\eqref{eq:COUPLING}.
Observe that this subspace can be represented according to
\[
  W_{\q} = w_{\q}+V_{\q} \subset H_{\q},
\]
where $w_{\q} \in H_{\q}$ takes care of prescribed contour and slope by satisfying  
\begin{equation} \label{eq:AFINEW}
w_{\q}(x) = h_i(\Phi_{\q_i}^{-1}(x)), \qquad {\textstyle  \frac{\partial}{\partial n}} w_{\q}(x) = s_i(x),\qquad x \in \Gamma_{\q_i},\quad i=1,\dots,N.
\end{equation}
In order to incorporate the remaining degrees of freedom representing height and tilt of particles, 
we use the direct sum
\[
  V_{\q} = H^2_0(\Omega_{\q}) + \spn \{ \eta_{i}^0, \eta_{i}^1, \eta_{i}^2 \;|\; i= 1,\dots, N \}
\]
of  $H^2_0(\Omega_{\q})=
\{v \in H^2(\Omega_{\q})\;|\; v = \frac{\partial}{\partial n} v = 0 \text{ on }\partial \Omega_\q \}$
and the $3N$-dimensional subspace spanned by functions $\eta_{i}^0, \eta_{i}^1, \eta_{i}^2  \in H_{\q}$ with the properties
\[
\eta_{i}^0(x) = 1, \quad \eta_{i}^1(x) = x_1 -X_{i,1},\quad \eta_{i}^2(x) = x_2 -X_{i,2},  \qquad x \in \Gamma_{\q_i},
\]
and $\eta_{i}^k = \frac{\partial}{\partial n}\eta_{i}^k = 0$ on $\Gamma_{\q_j}$ with $k=0,1,2$ and $i\neq j =1,\dots, N$.
Observe that  the additional degrees of freedom $\bgamma \in \R^{3N}$ are incorporated into the solution space in this way,
since all particle motions of the form \eqref{eq:height_tilt_function},
involving height $Z_i$ and tilt $\beta_i$, are precisely the linear combinations
\begin{align*}
  Z_i \eta_i^0 + \beta_{i,1} \eta_i^1 + \beta_{i,2}\eta_i^2.
\end{align*}

As $V_{\q}$ is a closed subspace of the Hilbert space $H_{\q}$ 
and the derivative of $\cJ_{\Omega_{\q}}$ is linear and $H_{\q}$-elliptic, 
the following existence and uniqueness result follows from the Lax--Milgram lemma.
\begin{theorem} \label{theo:MEMMIN}
Assume that $h_i \in H^{3/2}(\Gamma_i)$ and  $s_i \in H^{1/2}(\Gamma_i)$, $i=1,\dots, N$, and $\q \in \cD$. 
Then  there is $w_{\q}\in H_{\q}$ with the property \eqref{eq:AFINEW}
and the minimization problem
\begin{equation}  \label{eq:MEMMIN}
 u_{\q} = \argmin{v \in W_{\q} } \cJ_{\Omega_{\q}}(v)
\end{equation}
  has a unique solution $u_{\q} \in W_{\q}\subset H_{\q} \subset H^2(\Omega_{\q})$.
\end{theorem} 
It is well-known that \eqref{eq:MEMMIN} is equivalent to the variational equation
\begin{equation}\label{eq:MEMVAR}
  u_{\q} \in  w_{\q} + V_{\q}:\qquad
  \int_{\Omega_{\q}}\kappa \Delta u_{\q}\Delta v + \sigma \nabla u_{\q} \cdot  \nabla v\; dx = 0
\qquad \forall v \in V_{\q}.
\end{equation}
which in turn can be regarded as a parametrized linear elliptic partial differential equation of fourth order.
Note that contributions from $H^2_0(\Omega_{\q})$ and  $W_{\q}$ to the solution $u_{\q}$ 
can be decoupled by orthogonalization, cf. \cite[Lemma~3.1]{GraeserKies2017a}.
This turns out to be beneficial both for further analysis and  finite element approximation.

With $u_{\q}$ denoting the unique solution of \eqref{eq:MEMMIN} for given $\q \in \cD$,
Theorem~\ref{theo:MEMMIN} now allows to define the {\em geometric interaction potential} of particles in membranes 
according to
\begin{equation} \label{eq:GEOINT}
\cM(\q) = \cJ_{\Omega_{\q}}(u_{\q}) = \underset{ v \in  W_{\q} }{\min}\cJ_{\Omega_{\q}}(v),\qquad \q \in \cD.
\end{equation}
For the construction and analysis of finite element approximations of the minimization problem \eqref{eq:MEMMIN}, 
and thus of $\cM(\q)$, we refer, e.g., to~\cite{ElliottEtAl2016,GraeserKies2017a,GraeserKies2019,KiesT19}.

\subsection{Differentiability and stable representation of gradient}
It has been shown by Kies~\cite{KiesT19} and Gr\"aser and Kies~\cite[Lemma~4.6]{GraeserKies2017a} 
that the  geometric interaction potential $\cM: \cD \subset \R^{3 N}\mapsto \R$ defined in \eqref{eq:GEOINT} is a smooth function.
\begin{theorem}
The geometric interaction potential $\cM$ defined in \eqref{eq:GEOINT}  is differentiable in a neighborhood of $\q \in \cD$.
\end{theorem}
In particular, this implies existence of the gradient $\nabla \cM(\q)$,
\[
\nabla \cM(\q) = \left(  \partial_{\q_i} \cM(\q) \right)_{i=1}^N, 
\quad  \partial_{\q_i} \cM(\q) = ( \partial_{X_i} , \partial_{\alpha_i})\cM(\q),
\qquad \q \in \cD.
\]
We are interested in a computationally feasible representation of $\nabla \cM(\q)$. To this end, we consider
the directional derivative $\partial_{\ce}\cM(\q)$ at $\q \in \cD$ in an arbitrary direction $\ce \in \R^{3 N}$.
Such a representation was recently  provided by Kies~\cite{KiesT19} and Gr\"aser and Kies~\cite[Lemma~4.8]{GraeserKies2017a}
utilizing the framework of shape calculus~\cite{sokolowskizolesio92}.
For each fixed direction 
\[\
\ce = (\ce_i)_{i=1}^N\in \R^{3N}, \quad \ce_i=(E_i, \delta_i) \in \R^2\times \R,
\]
it requires the construction of a twice differentiable vector field $\phi_{\ce}: \Omega_{\q} \mapsto \R^2$ with the property
\begin{equation} \label{eq:VPROP}
  \phi_{\ce}(x) = E_i + \delta_i Q \left (x - X_i \right), \quad D \phi_{\ce}(x) = \delta_i Q, \qquad x \in \Gamma_{\q_i},
\end{equation}
with  $Q =  \left(  \begin{matrix} 0 & -1\\ 1 & 0   \end{matrix}\right)$ and $D \phi_{\ce}$ denoting the Jacobian matrix of $\phi_{\ce}$.
Notice that this property ensures that the restriction
of the vector field $\phi_{\ce}$ to $\Gamma_{\q_i}$ coincides with the
directional derivative of the particle rigid body motion $\Phi_{\q_i}(y)$ at $\q_i$ in the direction $\ce_i$.
Utilizing the notation $A:B$ for the Frobenius inner product
of two matrices $A$, $B\in \R^{2\times 2}$,
$I \in \R^{2 \times 2}$ for the identity matrix in $\R^2$,
and  $D^2 u_{\q}$   for the Hessian of $u_{\q}$, 
the desired representation of the directional derivative $\partial_{\ce} \cM(\q)$ reads as follows.
\begin{theorem}
Let $\q \in \cD$,  $u_{\q}$  be the solution of \eqref{eq:MEMMIN},  
$\phi_{\ce}: \Omega_{\q} \mapsto \R^2$ be a vector field with property \eqref{eq:VPROP},
and $\cA =  \operatorname{div} (\phi_{\ce}) I - D \phi_{\ce} - D \phi_{\ce}^T$.
Then 
\begin{equation} \label{eq:GRADREP}
\begin{array}{c}
\displaystyle \partial_{\ce} \cM(\q) = 
  \int_{\Omega_{\q}} \kappa \Delta u_{\q} \left( \cA : D^2u_{\q} - \Delta \phi_{\ce} \cdot \nabla u_{\q} - \half \operatorname{div} (\phi_{\ce})\Delta u_{\q} \right) \; dx  \\
\displaystyle  + \half  \int_{\Omega_{\q}} \sigma \cA \nabla u_{\q} \cdot \nabla u_{\q} \; dx .
\end{array}
\end{equation}
\end{theorem}

This representation  \eqref{eq:GRADREP} does not require any additional regularity of the solution  $u_{\q}$ of \eqref{eq:MEMMIN}.
Moreover,  $\partial_{\ce}\cM(\q)$ depends continuously on  $u_{\q}$ with respect to the $H^2$-norm.
As a consequence, discretization error estimates for suitable finite element approximations of $u_{\q}$ directly
carry over to $\partial_{\ce}\cM(\q)$. We refer to \cite{GraeserKies2019,KiesT19} for details.
Such kind of properties are not available for straightforward finite difference approximations.

Suitable vector fields $\phi_{\ce}$ can be easily constructed, e.g., by one of the following two algorithms.
For simplicity, we fix
$i=1,\dots,N$ and let  $\ce_j = 0$ for $i\neq j$
so that $\partial_{\ce}\cM(\q)$ becomes the partial derivative
of $\cM$ with respect to a single particle motion.

The first algorithm starts with the selection of a closed neighborhood
of $\Gamma_{\q_i}$ of thickness $\varepsilon>0$ that does neither
intersect any other $\Gamma_{\q_j}$ nor the boundary $\partial \Omega$.
Now let $\xi: \Omega_{\q} \to [0,1]$ be a twice differentiable scalar function with
\begin{align*}
  \xi(x) 
  \begin{cases}
    = 1 & \text{if } \operatorname{dist}(x,\Gamma_{\q_i})< \varepsilon/2,\\
    =0 & \text{if } \operatorname{dist}(x,\Gamma_{\q_i})> \varepsilon,\\
    \in [0,1] & \text{else}.
  \end{cases}
\end{align*}
Then the vector field
\begin{align*}
  \phi_{\ce}(x) = \xi(x) \Bigl(E_i + \delta_i Q \left (x - X_i \right) \Bigr)
\end{align*}
has the desired properties. Notice that the condition on the
Jacobian matrix $D \phi_{\ce}$ is satisfied, since $\phi_{\ce}$
is affine linear in the $\varepsilon/2$ neighborhood of $\Gamma_{\q_i}$.
The complexity of implementing this algorithm depends on the
accessibility of the particle shape and a suitable $\varepsilon$-neighborhood.
For example, in case of circular particles $B_i$, $\varepsilon$ can be computed
from the particle diameters and distances and $\xi(x)$ can be chosen rotationally
symmetric.

As an alternative that is suitable for complex particle shapes, one could compute a vector field
$\phi_{\ce}$ by solving the partial differential equation
\begin{align*}
  (\kappa \Delta^2 - \sigma \Delta) \phi_{\ce}(x) &= 0 & x &\in \Omega_{\q},\\
  \phi_{\ce}(x) = \frac{\partial}{\partial n} \phi_{\ce}&= 0 & x &\in \partial\Omega
\end{align*}
with the additional boundary condition~\eqref{eq:VPROP} on the
particle boundaries.
Up to the boundary conditions, this equation is of the same
nature as~\eqref{eq:MEMVAR} and can be discretized using
the same finite element techniques.
While finite element approximations will typically only be
weakly but not strongly twice differentiable, this is still
feasible, since the expression~\eqref{eq:GRADREP} is
$H^2$-continuous also with respect to $\phi_{\ce}$.

\subsection{Direct interaction potential}
The membrane-mediated geometric interaction potential is augmented by direct \emph{particle--particle interaction}
as expressed by a potential $\cP(\q)$. 
There is a wide variety of such potentials depending on the properties of the particles under consideration.
Here, we only consider so-called {\em soft-wall constraints}
that could be used to incorporate the condition $\q \in \cD$ by penalization~\cite{ElliottEtAl2016}
and take the form
\begin{equation}\label{eq:PARTPART}
 \cP=\cP_1+\cP_2.
\end{equation}
The first contribution $\cP_1$ consists of a Lennard-Jones-type potential 
\[ 
    \cP_1(\q)=\sum_{\underset{i\neq j}{i,j=1}}^N \cP_{ij},\;\;
    \cP_{ij}=4 \epsilon_{ij}\left[ \left(\frac{\sigma_{ij}}{\text{dist}(B_{\q_i},B_{\q_j})}\right)^{12}-\left(\frac{\sigma_{ij}}{\text{dist}(B_{\q_i},B_{\q_j})}\right)^{6} \right],
\]
such that $\text{dist}(B_{\q_i},B_{\q_j}) > 0$ for $i\neq j$ and $\cP_1(\q)=\infty$ otherwise.
This term  accounts for  the repulsion and attraction of particles. Similarly, we set 
\[
    \cP_2(\q)=\sum_{i=1}^N \left(\frac{\sigma_{i}}{\text{dist}(B_{\q_i},\partial \Omega)}\right)^6
\]
such that $\text{dist}(B_{\q_i},\partial \Omega)>0$, $i=1,\dots,N$, and $\cP_2(\q)=\infty$ otherwise.  
This term is accounting for escaping particles. 
For circular particles $B_{\q_i}$ with radius $r_i$, we have $\text{dist}(B_{\q_i},B_{\q_j})=|X_i-X_j|-(r_i+r_j)$.
Note that the soft-wall potential $\cP=\cP_1 + \cP_2$ is continuously differentiable on $\cD$.

\subsection{Full interaction potential}
The full interaction potential of particles in membranes $\cE(\q)$ finally reads
\begin{equation} \label{eq:FULLPOT}
 \cE(\q)= \cM(\q) + \cP(\q),\qquad \q \in \cD.
\end{equation}
Membraned-mediated clustering of particles  is often  associated with local or global minima of 
the full particle potential $\cE$ on $\cD\subset \R^{3N}$ that can be computed, e.g., 
by gradient-related  optimization methods (cf., e.g. \cite{KiesT19,delleSiteGraeserKuscheKornhuber20}).

\section{Langevin sampling of Helmholtz free energy differences}
\subsection{Langevin dynamics}
As a starting point for particle dynamics, we introduce the \emph{separable Hamiltonian}
\begin{equation} \label{eq:HAMILTONIAN}
  H(\q,\vel) \colonequals \cE(\q) + \half \vel \cdot \vel
\end{equation}
of $N$  particles with locations $\q =(\q_i)_{i=1}^N \in \R^{3N}$,  velocities $\vel = ({\vel}_i)_{i=1}^N \in \R^{3N}$, 
and the interaction potential $\cE$ of particles introduced in  \eqref{eq:FULLPOT}.
Assuming that the membrane particle system is embedded in an infinite heat bath that keeps the temperature $\Te$ constant,
we consider the associated stochastic Langevin process~\cite{tuckbook}.

\begin{align}
  \begin{split}
    \label{eq:langevin}
    d{\q_t} & = \vel_t d{t}\\
    d{\vel}_t & = -\nabla \cE(\q_t)d{t}- \fric \vel_t d{t} + \sqrt{2 \gamma / \beta } d{\qB_t}
  \end{split}
\end{align}
denoting $\beta =  (k_B \Te)^{-1}$.
Thermal fluctuations of  particles are represented by scaled $3N$-dimensional Brownian motion $\qB_t$,
and the corresponding thermal energy  is balanced by a viscous friction term $\fric \vel$.
Observe that we consider thermal fluctuations only of the particles (within the membrane),
but not of the membrane itself (within the surrounding heat bath).
This can be justified by the assumption that the membrane is 'fast' in the sense that
both fluctuations and particle-induced deformations of the membrane happen on much smaller time scales than particle motion.

After rescaling time according to  $t \mapsto \gamma t$ and letting $\gamma$ tend to infinity,
we formally obtain the overdamped Langevin equation 
\begin{equation} \label{eq:OLA}
d{\q_t} = - \nabla \cE(\q_t)dt + \sqrt{2 / \beta } d{\qB_t}, 
\end{equation}
where the rescaled particle trajectories and Brownian motion are still denoted by $\q_t$ and $\qB_t$, respectively.
The system \eqref{eq:OLA} of stochastic differential equations is completed by initial conditions
\begin{equation} \label{eq:INIT}
\q_0 \in \cD.
\end{equation}
Modeling fluctuating particle positions $\q_i=(X_i, \alpha_i)$ by Brownian motion
is well-established as far as spatial  coordinates $X_i$ are concerned
but seems to be less common for rotation angles $\alpha_i$.

Discretization of \eqref{eq:OLA} on a given time interval $[0,T]$ with given $T>0$ 
can be performed by the  Euler--Maruyama scheme
\begin{align}
\label{eq:eulerMaruyama}
  \qA_{k+1}= \qA_{k}-\tau\nabla \cE(\qA_{k}) + \sqrt{2\tau / \beta} \qG_{k},\qquad k = 0, 1, \dots , M,
\end{align}
with uniform time step size $\tau=T/M$ for some $M\in \N$, 
and independent, identically distributed centered Gaussian random variables $\qG_{k} \in \R^{3N}$ with unit variance.

For the interaction energy $\cE = \cM + \cP$  defined in \eqref{eq:FULLPOT} 
with geometric interaction $\cM$ taken from \eqref{eq:GEOINT}
existence and uniqueness of a discrete solution $\qA_{k+1}$
is guaranteed by Theorem~\ref{theo:MEMMIN}, 
provided that $\qA_{k}\in \cD$. 
This property can be guaranteed in various ways, for example by imposing additional reflection conditions.
Note that each time step of the  Euler--Maruyama scheme \eqref{eq:eulerMaruyama} 
requires the approximate solution of a  partial differential equation of the form \eqref{eq:MEMVAR},
e.g., by finite elements~\cite{GraeserKies2019},
in order to approximately evaluate the gradient $\nabla \cM(\qA_{k})$ 
via the representation formula \eqref{eq:GRADREP}.

\subsection{Energy perturbation method} \label{sec:APPROXFED}
We consider the canonical ensemble of  a fixed number of  particles with Hamiltonian $H$
 in a heat bath with fixed temperature $\theta > 0$ and fixed volume.
In order to deduce macroscopic properties from a given microscopic observable $A$,
\[
\TD \ni (\q,\vel) \mapsto A(\q,\vel) \in \R,
\]
defined on phase space $\TD = \cD \times \R^{3N}$ of all feasible pairs $\left(\q,\vel\right) = \left((\q_i, \vel_i)\right)_{i=1}^N$ 
of locations $\q\in \cD$ and velocities $\vel\in \R^{3N}$,
we introduce the expectation
\[
  \EXP_\mu(A) = \int_{\TD} A(\q,\vel)\; d\mu(\q,\vel)
\]
with respect to the canonical measure
\begin{align*}
 d \mu(\q,\vel) = Z_\mu^{-1} e^{-\beta H(\q,\vel)} d(\q,\vel)  ,\qquad
  Z_\mu = \int_{\TD} e^{-\beta H(\q,\vel)} d(\q,\vel).
\end{align*}

As a related macroscopic quantity, we consider the absolute Helmholtz free energy defined (up to a constant) by
\begin{align*}
  F = -\frac{1}{\beta} \ln(Z_\mu).
\end{align*}
It can be interpreted as a measure of stability  of the system
in the sense that lower Helmholtz free energies are expressing more stable macroscopic states.

Assuming that the Hamiltonian $H = H_{\omega}$ is parametrized by some reaction coordinate $\omega\in [-1,1] \subset \R$,
we aim at the approximation of  Helmholtz free energy differences $\Delta F({\omega})  = F_{\omega} - F_{0}$, $\omega\in [-1,1]$.
Utilizing the free energy perturbation method due to Zwanzig~\cite{zwanzig1954high}, \cite[Section~1.3]{lelievre},
 $\Delta F({\omega})$ can be expressed as the expectation of the observable 
 \begin{align}\label{eq:zwanzig_observable}
  A(\q,\vel) = e^{-\beta(H_{\omega}(\q,\vel)-H_{0}(\q,\vel))}
 \end{align}
 according to
\begin{equation} \label{eq:HFEGEN}
\begin{array}{rl}
 \displaystyle  e^{-\beta \Delta F({\omega})}
    &  \displaystyle  = Z_{\mu_0}^{-1} \int_{\TD} e^{-\beta H_{\omega}(\q,\vel)} \;d{(\q,\vel)}\\
    &  \displaystyle = \int_{\TD}A(\q,\vel)\; d\mu_0(\q,\vel)
       = \EXP_{\mu_0}(A) 
    \end{array}
\end{equation}
with $\mu_0$ denoting the canonical measure induced by $H_{0}$.

We  assume that  $H_{\omega}$ takes the form
\[
H_{\omega}(\q,\vel) = \cE_{\omega}(\q) + \half \vel \cdot \vel, \qquad \omega\in [-1,1],
\]
or, more generally, is separable with  kinetic energy independent of $\omega$.
Then the observable $A$ from~\eqref{eq:zwanzig_observable} satisfies
\begin{equation} \label{eq:REDOBS}
A(\q,\vel) = e^{-\beta(\cE_{\omega}(\q)-\cE_{0}(\q))} = A(\q),
\end{equation}
i.e., it is independent of $\vel$ and the  representation \eqref{eq:HFEGEN} reduces  to
\begin{align} \label{eq:REDHFED}
  e^{-\beta\Delta F({\omega})}
  =  \EXP_{\mu_0}(A) 
  = Z_{\nu_0}^{-1} \int_{\cD} e^{-\beta (\cE_{\omega}(\q)-\cE_{0}(\q))} e^{-\beta\cE_0(\q)}\; d\q
    = \EXP_{\nu_0}(A)
\end{align}
with the reduced canonical measure $d \nu_0$ given by
\[
d \nu_0(\q) = Z_{\nu_0}^{-1} e^{-\beta \cE_{0}(\q)}\; d \q ,
\qquad Z_{\nu_0} = \int_{\cD} e^{-\beta \cE_{0}(\q)}\; d \q.
\]

\subsection{Langevin sampling} \label{subsec:SAMPLER}
We now assume ergodicity with respect to $d \nu_0$ in the sense that
\begin{align}   \label{eq:ergodic}
  \EXP_{\nu_0}(A) = \lim_{T\to\infty}  \;\frac{1}{T} \int_0^T A(\q_t)\; dt \qquad \text{a.s.}
\end{align}
holds  for every observable $A(\q)$ with  $\q_t$ satisfying the overdamped Langevin equation~\eqref{eq:OLA} 
with the reference energy $\cE = \cE_0$.  This essentially means that samples of the trajectory $\q_t$ are visiting 
the whole configuration space~$\cD$ while replicating the density of the canonical measure~$d \nu_0$.
Together with \eqref{eq:REDOBS} and \eqref{eq:REDHFED} ergodicity \eqref{eq:ergodic} leads to the representation
\begin{equation} \label{eq:FEDIR}
e^{-\beta\Delta F({\omega})} =  \lim_{T\to\infty} \; \frac{1}{T} \int_0^T  e^{-\beta(\cE_{\omega}(\q_t)-\cE_{0}(\q_t))}  \; dt.
\end{equation}
Formally approximating  the integral by a Riemannian sum with stepsize $\tau = T/M$ and the stochastic process $\q_t$, $t \in [0,T]$,
by a realization $\qA = (\qA_{k})_{k=0}^M$ of  the Euler--Maruyama discretization~\eqref{eq:eulerMaruyama}
for $\cE = \cE_{0}$ and the same stepsize $\tau$,
we obtain the  Langevin sampler
\begin{equation} \label{eq:LANGEVIN}
  \widehat{A}_M =  \frac{1}{M+1} \sum_{k=0}^M   e^{-\beta\left(\cE_{\omega}( \qA_{k})  - \cE_{0}(\qA_{k})\right)}
\end{equation}
for the expectation value $\EXP_{\nu_0}(A)$.
Together with \eqref{eq:REDHFED} this provides the sampler
\begin{equation} \label{eq:LASA1}
  \widehat{\Delta F}_M({\omega}) =  - \frac{1}{\beta} \ln (\widehat{A}_M)
\end{equation}
of the free energy difference $\Delta F({\omega})$.

Following  \cite[Section 2.3.1]{lelievre}, the sampling error of $\widehat{A}_M$
can be decomposed into the bias and the statistical error according to
\begin{align*}
  \EXP(|\widehat{A}_M - \EXP_{\nu_0}(A)|^2)
    = \Big(\EXP(\widehat{A}_M) - \EXP_{\nu_0}(A)\Big)^2
      + \EXP\Big(|\widehat{A}_M - \EXP(\widehat{A}_M)|^2\Big).
\end{align*}
Assuming that the discrete stochastic process $(\qA)_{k=0}^\infty$ 
obtained from the Euler--Maruyama scheme~\eqref{eq:eulerMaruyama}
samples some invariant measure $\tilde{\nu}_0$, the bias  can be estimated by
\begin{align*}
    \Big|\EXP(\widehat{A}_M) - \EXP_{\nu_0}(A)\Big|
    \leq
    \Big|\EXP(\widehat{A}_M) - \EXP_{\tilde{\nu}_0}(A)\Big|
    + \Big|\EXP_{\tilde{\nu}_0}(A) - \EXP_{\nu_0}(A)\Big|.
\end{align*}
Under the  assumption that $\EXP(A(\qA_k))$ converges exponentially fast to
$\EXP_{\tilde{\nu}_0}(A)$ as $k\to \infty$, one can show that
the finite sampling bias $|\EXP(\widehat{A}_M) - \EXP_{\tilde{\nu}_0}(A)|$ is of order $\cO(M^{-1})$.
The  second term, the perfect sampling bias $|\EXP_{\tilde{\nu}_0}(A) - \EXP_{\nu_0}(A)|$,
is associated with the error of time discretization, which, under  suitable assumptions,
is of  order $\cO(\tau)$ for the actual Euler--Maruyama scheme 
and of higher order for more advanced time discretizations (see, e.g., \cite{abdulle2014high} and the references cited therein).

The remaining statistical error $\EXP\Big(|\widehat{A}_M - \EXP(\widehat{A}_M)|^2\Big)$ typically satisfies a central limit theorem of the form
\begin{align*}
  \sqrt{M}|\widehat{A}_M - \EXP(\widehat{A}_M)| \to \mathcal{N}(0,\sigmaLangevin^2),
\end{align*} 
with variance $\sigmaLangevin$ associated with the actual discrete Langevin sampling strategy. 
Note that $\sigmaLangevin$ is expected to be much smaller than, e.g., the variance of standard Monte-Carlo sampling.

\section*{Acknowledgement} This research has been funded by Deutsche Forschungsgemeinschaft (DFG) through the grant CRC 1114: ``Scaling Cascades in Complex Systems'', Project Number 235221301, Project A07.

\bibliographystyle{elsarticle-num} 
\bibliography{DSGKKarxive1}





\end{document}